\documentclass[12pt, reqno]{amsart}
\usepackage[utf8]{inputenc}
\usepackage{url}
\usepackage{amsmath,color}
\usepackage{amsfonts}
\usepackage{amssymb}
\usepackage{amsthm}
\usepackage{multirow}
\usepackage{graphicx} 
\usepackage{lipsum}
\usepackage{xcolor}

\newtheorem{theorem}{Theorem}
\newtheorem{remark}[theorem]{Remark}
\newtheorem{lemma}[theorem]{Lemma}

\newtheorem{proposition}[theorem]{Proposition}

\newcommand{\tmop}[1]{\ensuremath{\operatorname{#1}}}
\numberwithin{equation}{section}
\numberwithin{theorem}{section}

\numberwithin{equation}{section}

\def\sumstar{\sideset{}{^*}\sum}

\def\pamod{\!\!\!\!\pmod}

\newcommand{\fixmehide}[1]{}
\newcommand{\fixmelater}[1]{}
\newcommand{\fixmedone}[1]{}
\newcommand{\fixmehidden}[1]{}

\newcommand{\eps}{\varepsilon}

\title{The signs of products of central $L$-values of twists of  modular forms}

\author{Hung M. Bui, Natalie Evans, Stephen Lester, and Kyle Pratt}
\address{Department of Mathematics, University of Manchester, Manchester M13 9PL, UK}
\email{hung.bui@manchester.ac.uk}
\address{Department of Mathematics, King’s College London, London WC2R 2LS, UK}
\email{natalie.r.evans@kcl.ac.uk}
\address{Department of Mathematics, King’s College London, London WC2R 2LS, UK}
\email{steve.lester@kcl.ac.uk}
\address{Brigham Young University, Department of Mathematics, Provo, UT 84602, USA}
\email{kyle.pratt@mathematics.byu.edu}

\allowdisplaybreaks

\begin{document}

\begin{abstract}
Let $f,g$ be distinct weight $2$ newforms of fixed level whose product of root numbers equals $1$ and $\chi$ be a Dirichlet character modulo a prime $q$. We study the sign of the product of the central $L$-values of $f$ twisted by $\chi$ and $g$ twisted by $\overline \chi$. We prove that the number of twists for which the product of the two central $L$-values is positive tends to infinity with $q$ and the number of twists for which it is negative also grows with $q$. Moreover, assuming the Generalized Riemann Hypothesis we obtain a positive proportion of twists for which the product is positive and a positive proportion of twists for which it is negative. 

\end{abstract}

\maketitle

\section{Introduction}

Central values of $L$-functions are fundamental objects of study in number theory, partly due to the connection to the Birch and Swinnerton-Dyer Conjecture. It is widely expected that the vanishing of an $L$-function at the center of the critical strip, apart from trivial reasons such as the functional equation having sign $-1$, reflects deep arithmetic properties of the underlying object to which the $L$-function is attached.

Let $f$ be a weight $2$ newform on $\Gamma_0(N)$ and $q$ be a prime with $(q,N)=1$. For $\chi$ a primitive character modulo $q$ and $\tmop{Re}(s) >1$, the $L$-function associated to $f\otimes \chi$ is given by
\[
L(s, f \otimes \chi) =\sum_{n \ge 1} \frac{\lambda_f(n)\chi(n)}{n^s}=\prod_p\bigg(1-\frac{\alpha_{f,1}(p) \chi(p)}{p^s} \bigg)^{-1}\bigg(1-\frac{\alpha_{f,2}(p)\chi(p)}{p^s} \bigg)^{-1}.
\]
Here and throughout $\lambda_f(n)$ is the $n$th Hecke eigenvalue of $f$ and $\alpha_{f,j}$, with $j=1,2$, are the Satake parameters, which satisfy
\[
\alpha_{f,1}(p)+\alpha_{f,2}(p)=\lambda_f(p) \qquad \text{and} \qquad \alpha_{f,1}(p)\alpha_{f,2}(p)=\chi_N(p),
\]
where $\chi_N$ is the principal character modulo $N$.  Building on works of Rohrlich \cite{R1,R2}, Chinta \cite{Ch} showed that $L(\frac12,f\otimes \chi) \neq 0$ for almost all primitive characters $\chi \pmod q$. Chinta's result is, however, insufficient to rule out the possibility of extremely small central values. Using the mollification method, Blomer et al. \cite[Theorem 1.8]{BFK} established that
\[
\#\Big\{\chi \pamod q\ \text{primitive}\,:|L(\tfrac12,f\otimes\chi)|\geq \frac{1}{\log q}\Big\}\gg q.
\]
Under the assumption of the Generalized Riemann Hypothesis (GRH), the authors \cite[Corollary 1.5]{BELP} also obtained
positive proportions of twists for which the central $L$-values associated to two distinct modular forms simultaneously grow or shrink with $q$. In this work, we are interested in the signs of the products of those two central $L$-values.

\subsection{Main results}  Throughout the paper, we let $f,g$ be two distinct weight $2$ newforms on $\Gamma_0(N)$ and assume that $N$ is fixed. We write $\eps(f) \in \{ \pm 1\}$ for the root number of $L(s,f)$, see \eqref{eq:ffunctionaleq}. As a consequence of the functional equation for $L(s,f\otimes \chi)$ stated in \eqref{eq:fchifunctionaleq}, see also \eqref{eq:athalf}, if $\eps(f)\eps(g)=1$ then
\begin{equation} \label{eq:real}
L(\tfrac12,f\otimes\chi)L(\tfrac12,g\otimes\overline{\chi})\in\mathbb{R}.
\end{equation}
Since for complex $\chi$ and real $s$, $L(s,f \otimes \chi) \overline{L(s,g \otimes \chi)}$ is not real-valued, a negative central value does not yield a counterexample to GRH, whereas a negative central $L$-value in the case where $\chi$ is the quadratic character modulo $q$ would. 
 In this article we show that for $f,g$ with $\eps(f)\eps(g)=1$, the number of twists $\chi \pmod q$ for which the product of central $L$-values is positive tends to infinity with $q$ and the number of twists for which it is negative also grows with $q$. Moreover, under GRH we obtain positive proportions of such twists. 

\begin{theorem} \label{thm0}
Suppose $\eps(f)\eps(g)=1$. Then for any $\delta>0$ we have
\[
\#\Big\{\chi \pamod q\ \tmop{primitive}\,: L(\tfrac12,f\otimes\chi)L(\tfrac12,g\otimes\overline{\chi})>0\Big\} \gg_{\delta} \, q^{1/4-\delta}
\]
and 
\[
\#\Big\{\chi \pamod q\ \tmop{primitive}\,: L(\tfrac12,f\otimes\chi)L(\tfrac12,g\otimes\overline{\chi})<0\Big\} \gg_{\delta} \, q^{1/4-\delta}.
\]
\end{theorem}
\begin{remark}
    \emph{In \cite{BFK}, Blomer et al. compute the first moment of the product of the two automorphic $L$-functions. For $\eps(f)\eps(g)=1$ the main term in the first moment \cite[Proposition 5.2]{BFK} is $\asymp q$ and is positive. Hence, the first claim of Theorem \ref{thm0} follows from \cite[Proposition 5.2]{BFK} along with the subconvexity estimate
    \begin{equation} \label{eq:subconvexbd}
    L(\tfrac12,f \otimes \chi) \ll q^{3/8+\eps}
    \end{equation}
    due to Blomer and Harcos \cite[Theorem 2]{BH}. However, it is not possible to deduce the second claim of Theorem \ref{thm0} by this argument.
    }
\end{remark}

\begin{theorem}\label{thm1}
Assume GRH and suppose $\eps(f)\eps(g)=1$. We have 
\[
\#\Big\{\chi \pamod q\ \tmop{primitive}\,: L(\tfrac12,f\otimes\chi)L(\tfrac12,g\otimes\overline{\chi})>0\Big\} \gg q
\]
and 
\[
\#\Big\{\chi \pamod q\ \tmop{primitive}\,: L(\tfrac12,f\otimes\chi)L(\tfrac12,g\otimes\overline{\chi})<0\Big\} \gg q.
\]
\end{theorem}

\begin{remark}
    \emph{In \cite{GI}, Gonek and Ivi\'c unconditionally established a similar result for the signs of the Hardy $Z$-function $Z(t)$. They used Selberg's mollifier, which does not work in our case, see the discussion in Section \ref{discussion}.}
\end{remark}

Theorem \ref{thm1} can be strengthened as follows.
\begin{theorem}\label{thm2}
Assume GRH and suppose $\eps(f)\eps(g)=1$. There exists $c>0$, which depends at most on $N$, such that
\[
\# \Big\{\chi \pamod q\ \emph{primitive}\,:L(\tfrac12,f\otimes\chi)L(\tfrac12,g\otimes\overline{\chi})>\exp\big(-c\sqrt{\log\log q}\big)\Big\}\gg_c q
\]
and
\[
\# \Big\{\chi \pamod q\ \emph{primitive}\,:L(\tfrac12,f\otimes\chi)L(\tfrac12,g\otimes\overline{\chi})<-\exp\big(-c\sqrt{\log\log q}\big)\Big\}\gg_c q.
\]
\end{theorem}

\begin{remark}
    \emph{Blomer et al. \cite[Theorem 1.12]{BFK} used Soundararajan's resonance method 
to unconditionally show the existence of twists $\chi \pmod q$ such that $|L(\tfrac12,f\otimes\chi)L(\tfrac12,g\otimes\overline{\chi})|$ is as large as $\exp\big(c\sqrt{\frac{\log q}{\log\log q}}\big)$. It is, however, not possible to remove the modulus sign in their arguments.}
\end{remark}

In the complementary case where $\eps(f)\eps(g)=-1$, it follows from the functional equation, see \eqref{eq:athalf}, that\footnote{For the quadratic twist, one of the two $L$-functions in the product will have root number equal to $-1$ and will vanish.}
\begin{equation} \label{eq:imag}
iL(\tfrac12,f\otimes\chi)L(\tfrac12,g\otimes\overline{\chi})\in \mathbb{R}.
\end{equation}
By a straightforward modification of our arguments, if $\eps(f)\eps(g)=-1$, then the analogues of Theorems \ref{thm0}, \ref{thm1}, and \ref{thm2} hold for both positive and negative values of  $i L(\tfrac12,f\otimes\chi)L(\tfrac12,g\otimes\overline{\chi})$. In particular, we can show under GRH that if $\eps(f)\eps(g)=-1$, then
\[
\# \Big\{\chi \pamod q\ \text{primitive}\,: i L(\tfrac12,f\otimes\chi)L(\tfrac12,g\otimes\overline{\chi}) > 0\Big\} \gg q.
\]

\subsection{Main propositions} \label{sec:molandprop}
For $\tmop{Re}(s)>1$, the $L$-function associated to $f$ is
\[
L(s,f)=\sum_{n \ge 1} \frac{\lambda_f(n)}{n^s}=\prod_{p}\bigg(1-\frac{\alpha_{f,1}(p)}{p^s} \bigg)^{-1}\bigg(1-\frac{\alpha_{f,2}(p)}{p^s} \bigg)^{-1}.
\]
The function $L(s,f)$ satisfies the functional equation
\begin{equation} \label{eq:ffunctionaleq}
\Lambda(s,f)=\varepsilon(f) \Lambda(1-s,f),
\end{equation}
where $\eps(f)\in \{\pm1\}$ is the root number and
\[
\Lambda(s,f)=N^{s/2}L_{\infty}(s,f)L(s,f), \qquad L_{\infty}(s,f)=(2\pi)^{-s}\Gamma(s+\tfrac{1}{2}).
\]
For non-principal $\chi \pmod{q}$, the functional equation for $L(s,f\otimes \chi)$ is given by 
\begin{equation} \label{eq:fchifunctionaleq}
\Lambda(s,f\otimes \chi)=\eps(f\otimes \chi)\Lambda(1-s, f\otimes \overline \chi),
\end{equation}
where
\[
\Lambda(s,f\otimes \chi)=(q^2N)^{s/2} L_{\infty}(s,f) L(s,f\otimes \chi)
\]
and
\[
\eps(f\otimes \chi)=\eps(f) \chi(N) \tau_{\chi}^2, \qquad \tau_{\chi}=\frac{1}{\sqrt{q}} \sum_{a \pamod{q}} \chi(a) e\Big(\frac aq\Big), \qquad e(x)=e^{2\pi i x},
\]
see \cite[Lemma 2.1]{BFK}.
Hence, it follows that 
\begin{equation} \label{eq:athalf}
L(\tfrac12,f\otimes \chi)L(\tfrac12,g\otimes \overline \chi)= \eps(f)\eps(g) L(\tfrac12,f \otimes \overline \chi)L(\tfrac12, g\otimes \chi),
\end{equation}
from which \eqref{eq:real} and \eqref{eq:imag} follow.

Our arguments rely upon a mollifier that behaves like an Euler product, which we define as follows. Let $\eta>0$ be fixed and sufficiently small. For $0 \le j \le J$, let $\theta_j =\frac{\eta e^j}{(\log \log q)^5}$ and $\ell_j=2 \lfloor \theta_j^{-3/4}\rfloor$, where $J$ is chosen so that $\eta \le \theta_J \le e \eta$. This implies that $J \asymp \log \log \log q$. Additionally, set $y=q^{\theta_0}$ and $x=q^{\theta_J}$. Let $c_0>0$ be fixed and sufficiently large and define $I_0=(c_0,y]$. For each $j=1,\ldots, J$, let $I_j=(q^{\theta_{j-1}},q^{\theta_j}]$. For $j=0,1,\ldots,J$, define the completely multiplicative functions $w_j(n)$ and $a_{f,j}(n)$ by
\begin{equation*}
\begin{split}
w_j(p)=p^{-\frac{1}{\theta_j \log q}} \Big(1-\frac{\log p}{\theta_j \log q} \Big), \qquad a_{f,j}(p)&=\lambda_f(p) w_j(p).
\end{split}
\end{equation*}
Let, for $j=0,1,\ldots, J$,
\begin{equation*}
M_{\frac12,f,j}(\chi)=\sum_{\substack{p|n \Rightarrow p \in I_j \\ \Omega(n) \le \ell_j}}\frac{  (-1/2)^{\Omega(n)}\nu(n)a_{f,J}(n) \chi(n)}{\sqrt{n}},
\end{equation*}
where $\Omega(n)=\sum_{p^a || n} a$ and $\nu$ is the multiplicative function with $\nu(p^a)=\frac{1}{a!}$. Also, let
\begin{equation} \label{eq:molldef}
M_{\frac12,f}(\chi)=\prod_{j=0}^{J} M_{\frac12,f,j}(\chi).
\end{equation}
The function $M_{\frac12,f}(\chi)$ is constructed to mollify $L(\frac12,f\otimes \chi)^{1/2}$. Since $\sum_{j=0}^J \theta_j \ell_j \le 20 \eta^{1/4}$ we also have that
\begin{equation} \label{eq:mollifierbd}
|M_{\frac12,f}(\chi)|= \prod_{j=0}^{J} |M_{\frac12,f,j}(\chi)| \le q^{\sum_{j=0}^J \theta_j \ell_j} \le q^{20\eta^{1/4}}.
\end{equation}
Additionally, we write $\sumstar$ to denote a sum over primitive characters.


We shall deduce Theorem \ref{thm0} from the following propositions.

\begin{proposition}\label{prop1}
We have
\[
\sumstar_{\chi \pamod q} L(\tfrac12,f\otimes\chi)L(\tfrac12,g\otimes\overline{\chi})M_{\frac12,f}(\chi)^2M_{\frac12,g}(\overline{\chi})^2\asymp q.
\]
\end{proposition}

\begin{proposition}\label{prop2}
We have
\[
\sumstar_{\chi \pamod q} L(\tfrac12,f\otimes\chi)L(\tfrac12,g\otimes\overline{\chi})\big|M_{\frac12,f}(\chi)M_{\frac12,g}(\chi)\big|^2\ll q(\log q)^{-1/2+\varepsilon}.
\]
\end{proposition}


For the proof of Theorem \ref{thm1} we also need the following proposition.

\begin{proposition}\label{prop3}
Assume GRH. Let $k>0$ and suppose that $ \eta \lceil k \rceil$ is sufficiently small. We have
\[
\sumstar_{\chi\pamod  q} \big|L(\tfrac12,f\otimes\chi)M_{\frac12,f}(\chi)^2\big|^{2k}\ll_k q.
\]
\end{proposition}
More precisely, in Proposition \ref{prop3} we assume that GRH holds for $L(s,\chi)$, $L(s,f\otimes\chi)$, $L(s,g\otimes \chi)$, $L(s,\tmop{sym}^2 f \otimes \chi)$, and $L(s,\tmop{sym}^2 g \otimes \chi)$ for every character $\chi \pmod{q}$.

%

\subsection{Discussion of the proofs and the organization of the paper}\label{discussion}

Since $M_{\frac12,f}(\chi)$ should act like $L(\frac12, f\otimes \chi)^{-1/2}$, we expect that 
\[
L(\tfrac12,f\otimes \chi)L(\tfrac12,g\otimes \overline \chi)\big|M_{\frac12}(\chi)\big|^2
\]
approximates the signum function of $L(\tfrac12,f\otimes \chi)L(\tfrac12,g\otimes \overline \chi)$. Obtaining cancellation on average of the signum function would imply approximately half of the non-zero values are positive and the other half negative.

In Proposition \ref{prop2}, we show that the mean value $L(\tfrac12,f\otimes \chi)L(\tfrac12,g\otimes \overline \chi)|M_{\frac12}(\chi)|^2$ is small. (The presence of the mollifier is necessary here, as the mean value of $L(\tfrac12,f\otimes \chi) L(\tfrac12,g\otimes \overline \chi)$ is $\asymp 1$ for $f,g$ with $\eps(f)\eps(g)=1$ by \cite[Proposition 5.2]{BFK}.) From this, we would like to conclude that $L(\tfrac12,f\otimes \chi)L(\tfrac12,g\otimes \overline \chi)$ is frequently positive and also frequently negative. However, we also need to rule out that the small mean value does not arise due to the presence of many small terms in the sum. This is accomplished by Proposition \ref{prop1}, and, \eqref{eq:subconvexbd} unconditionally or Proposition \ref{prop3} under GRH. With these estimates in hand, we prove Theorems \ref{thm0} and \ref{thm1} in Section \ref{sec:thm1}.

In addition to Propositions \ref{prop1}, \ref{prop2}, and \ref{prop3}, to establish Theorem \ref{thm2} we need to also control how often the mollifier $|M_{\frac12}(\chi)|$ is large. Since the mollifier acts like an Euler product, it should have an approximate log-normal distribution as $\chi \pmod{q}$ varies, with mean zero and variance of order $\log \log q$. This suggests that $|M_{\frac12}(\chi)|> \exp(c \sqrt{\log \log q})$ for  $o(q)$ characters $\chi \pmod q$ if $c$ tends to infinity with $q$. In Lemma \ref{lem:mollifierbd} we establish this estimate and complete the proof of Theorem \ref{thm2} in Section \ref{sec:thm2}.

To establish Proposition \ref{prop3} we require a mollifier that is a Dirichlet polynomial that acts like an Euler product. This idea was first introduced in \cite{RS} and has been successfully used in various other contexts, see \cite{LR, DFL,GZ,BELP,BFN,J,BFM1,BFM}. The sharp bounds for all $k$ in Proposition \ref{prop3} are obtained by adapting the methods of Soundararajan \cite{S} and Harper \cite{H}. This is the only place where we use GRH. Consequently, applying the twisted first moment proved in \cite[Theorem 5.1]{BFK} to evaluate the first moments appearing in Proposition \ref{prop1} and \ref{prop2} leaves a complex expression. To resolve this difficulty, in Section \ref{sec:prop12} we use the methods of \cite{BELP} and put forward a random model for $L(s,f\otimes\chi)L(s,g\otimes \overline \chi)$. This allows us to analyze the complex expression arising from the twisted first moment in terms of Euler products that we compute term-by-term, with different estimates for primes in different ranges. We note that for the proof of Theorem \ref{thm0}, a simpler mollifier would suffice, which would allow us to streamline this part of the argument.

The proof of Proposition \ref{prop3} closely follows the methods of \cite{LR,BELP}, which are now well-established in the literature. In \cite[Proposition 1.3]{BELP} a result analogous to Proposition \ref{prop3} is proved for a related, yet different mollifier. Additionally, in \cite[Proposition 4.1]{LR} upper bounds for mollified moments for the family of quadratic twists of $f$ are proved with an analogue of the mollifier considered in this work. For these reasons, we will omit the proof of Proposition \ref{prop3}. \newline 


\textbf{AI Disclosure:}
No AI tools were used in this work.
\newline 

\textbf{Acknowledgments:} KP is partially supported by the National Science Foundation (DMS-2418328) and the Simons Foundation (MPS-TSM-00007959).

\section{Proof of Theorems \ref{thm0} and \ref{thm1}} \label{sec:thm1}

Using the first two main propositions we will now establish Theorem \ref{thm0}.

\begin{proof}[Proof of Theorem \ref{thm0} assuming Propositions \ref{prop1} and \ref{prop2}] Let 
\[
S_+:=\Big\{\chi \pamod q\ \text{primitive}\,: L(\tfrac12,f\otimes\chi)L(\tfrac12,g\otimes\overline{\chi})>0\Big\}
\]
and
\[
S_-:=\Big\{\chi \pamod q\ \text{primitive}\,:L(\tfrac12,f\otimes\chi)L(\tfrac12,g\otimes\overline{\chi})<0\Big\}.
\]
Additionally, let
\begin{equation} \label{eq:Ipmdef}
I_{\pm}=\sumstar_{\chi\in S_{\pm}} L(\tfrac12,f\otimes\chi)L(\tfrac12,g\otimes\overline{\chi})\big|M_{\frac12,f}(\chi)M_{\frac12,g}(\chi)\big|^2.
\end{equation}
Then we have that
\begin{align*}
I_+-I_-&=\sumstar_{\chi\pamod q} |L(\tfrac12,f\otimes\chi)L(\tfrac12,g\otimes\overline{\chi})|\big|M_{\frac12,f}(\chi)M_{\frac12,g}(\chi)\big|^2\\
&\geq \bigg|\ \sumstar_{\chi\pamod q} L(\tfrac12,f\otimes\chi)L(\tfrac12,g\otimes\overline{\chi})M_{\frac12,f}(\chi)^2M_{\frac12,g}(\overline{\chi})^2\bigg|\gg q,
\end{align*}
by Proposition \ref{prop1}. Additionally, by  Proposition \ref{prop2} we have that
\begin{align*}
&I_++I_-=\sumstar_{\chi\pamod q} L(\tfrac12,f\otimes\chi)L(\tfrac12,g\otimes\overline{\chi})\big|M_{\frac12,f}(\chi)M_{\frac12,g}(\chi)\big|^2=o(q).
\end{align*}
Hence, it follows that
\begin{equation}\label{estI+-}
I_+\gg q\qquad\text{and}\qquad I_-\gg q.
\end{equation}

Applying \eqref{eq:subconvexbd} and \eqref{eq:mollifierbd} we also have that
\[
|I_{-}| \ll  q^{3/4+40\eta^{1/4}+\eps} |S_{-}|.
\]
Taking $\eta, \eps$ sufficiently small in terms of $\delta$ and using \eqref{estI+-} we get $|S_{-}| \gg q^{1/4-\delta}$. A similar argument, which we omit, also gives that $|S_+| \gg q^{1/4-\delta}$.
\end{proof}

Arguing similarly, and also using Proposition \ref{prop3}, we shall now quickly deduce Theorem \ref{thm1}

\begin{proof}[Proof of Theorem \ref{thm1} assuming Propositions \ref{prop1},  \ref{prop2}, and \ref{prop3}] 
Let $I_{\pm}$ be as defined in \eqref{eq:Ipmdef}.
Applying the Cauchy-Schwarz inequality twice, we also have that
\begin{align*}
|I_+| \le |S_+|^{1/2}\bigg(\ \sumstar_{\chi\pamod q} \big|L(\tfrac12,f\otimes\chi)M_{\frac12,f}(\chi)^2\big|^4\bigg)^{1/4}\bigg(\ \sumstar_{\chi\pamod q} \big|L(\tfrac12,g\otimes\chi)M_{\frac12,g}(\chi)^2\big|^4\bigg)^{1/4}.
\end{align*}
Using Proposition \ref{prop3} and \eqref{estI+-}, we deduce that $|S_+|\gg q$. By a similar argument, which we omit, we have that $|S_-| \gg q$.
\end{proof}

\section{Proof of Theorem \ref{thm2} assuming Propositions \ref{prop1}, \ref{prop2}, and \ref{prop3}} \label{sec:thm2}

To obtain lower bounds on the product of central $L$-values we require upper bounds on the mollifier. We first control how often the factors with $1 \le j \le J$ in the product defining the mollifier in \eqref{eq:molldef} can be large.
\begin{lemma} \label{lem:secondmomentbd}
We have 
\[
\sum_{\chi \pamod{q}} \prod_{j=1}^J \big|M_{\frac12,f,j}(\chi)\big|^2 \ll q(\log \log q)^{O(1)}.
\]

\end{lemma}

\begin{proof}
By construction,
there exists $\{ b_f(n)\}_{n \in \mathbb N} \subset \mathbb R$ such that
\begin{equation} \label{eq:Dirichlet}
\prod_{j=1}^J M_{\frac12,f,j}(\chi)=\sum_{ \substack{n \le q^{\gamma} \\ p|n \Rightarrow y < p \le x}}\frac{b_f(n) \chi(n)}{\sqrt{n}}
\end{equation}
where $\gamma=\sum_{j=1}^J \ell_j \theta_j \le 20 \eta^{1/4} <1$ and $|b_f(n)|\le C^{\Omega(n)}$ for some $C>0$, by Deligne's bound. Hence, 
\[
\begin{split}
\frac{1}{\varphi(q)}\sum_{\chi \pamod{q}} \prod_{j=1}^J \big|M_{\frac12,f,j}(\chi)\big|^2 =&\sum_{ \substack{n \le q^{\gamma} \\ p|n \Rightarrow y < p \le x}}\frac{b_f(n)^2 }{n}\\
\le & \sum_{ \substack{p|n \Rightarrow y < p \le x}}\frac{C^{2\Omega(n)} }{n} =\prod_{y < p \le x}\bigg(1+O\Big(\frac1p\Big) \bigg)\ll (\log \log q)^{O(1)},
\end{split}
\]
and we obtain the lemma.
\end{proof}

By construction the factor $M_{\frac12,f,0}(\chi)$ in the product defining the mollifier behaves like an Euler-product and the next lemma allows us to control how often it takes on large values. To state it, define
\[
P_f(\chi)=\sum_{p \in I_0}\frac{\lambda_f(p)w_J(p) \chi(p)}{\sqrt{p}}.
\]
\begin{lemma} \label{lem:devbd}
We have for $y^{V^2/9}<q$ that
\begin{equation} \notag
\#  \Big\{\chi \pamod q\ \tmop{primitive}\, : |P_f(\chi)| 
\ge V \sqrt{\tfrac12 \log \log y} \Big\} \ll q e^{-V^2/9}.
\end{equation}
\end{lemma}
\begin{proof} For $\ell,n \in \mathbb N$ with $n \ge 2$, we have that
\begin{equation} \label{eq:comb}
\# \{ (p_1,\ldots,p_{\ell}) : p_i \text{ is prime } \forall i \text{ and } p_1 \cdots p_{\ell}=n \}=\ell! \nu(n).
\end{equation}
Consequently, it follows that
\begin{equation} \label{eq:ppower}
P_f(\chi)^{\ell}=\ell! \sum_{\substack{p|n \Rightarrow p \in I_0 \\ \Omega(n)=\ell}} \frac{ \nu(n)a_{f,J}(n) \chi(n)}{\sqrt{n}}.
\end{equation}
Hence, for $y^{\ell}<q$ we have that
\begin{equation} \label{eq:momentbd}
\begin{split}
\frac{1}{\varphi(q)}\sum_{\chi \pamod q} |P(\chi)|^{2\ell} &=(\ell!)^2\sum_{\substack{p|n \Rightarrow p \in I_0 \\ \Omega(n)=\ell}} \frac{ \nu(n)^2a_{f,J}(n)^2 }{n} \\
&\le (\ell!)^2\sum_{\substack{p|n \Rightarrow p \in I_0 \\ \Omega(n)=\ell}} \frac{ \nu(n) a_{f,J}(n)^2 }{n}=\ell!\bigg(\sum_{p \in I_0}\frac{\lambda_f(p)^2w_J(p)^2 }{p}\bigg)^{\ell},
\end{split}
\end{equation}
where we used \eqref{eq:comb} in the last step.

Using Chebyshev's inequality, \eqref{eq:momentbd}, and Deligne's bound we conclude for $y^{\ell} < q$ that
\[
\begin{split}
\#  \Big\{\chi \pamod q\ \tmop{primitive}\, : |P_f(\chi)| 
\ge V \sqrt{\tfrac12 \log \log y} \Big\}
\le & \frac{1}{V^{2\ell} (\frac12 \log \log y)^{\ell}} \sum_{\chi \pamod q} |P_f(\chi)|^{2\ell} \\
\le & \frac{\ell! \varphi(q)}{V^{2\ell} (\frac12 \log \log y)^{\ell}} 
\bigg(4 \sum_{p \in I_0}\frac1p \bigg)^{\ell}.
\end{split}
\]
Hence, applying Stirling's formula we have for $y^{\ell} < q$ that
\[
\#  \Big\{\chi \pamod q\ \tmop{primitive}\, : |P_f(\chi)| 
\ge V \sqrt{\tfrac12 \log \log y} \Big\} \ll q \Big(\frac{9 \ell}{V^2 e} \Big)^{\ell}
\]
and taking $\ell=\lfloor V^2/9 \rfloor$ completes the proof.
\end{proof}

Using Lemmas \ref{lem:secondmomentbd} and \ref{lem:devbd}, we will obtain an upper bound on the frequency of large values of $|M_{\frac12}(\chi)|$, which is the last main ingredient needed for the proof of Theorem \ref{thm2}.
\begin{lemma} \label{lem:mollifierbd}
For $1\le \Lambda  \le (\log \log q)^{1/4}$ we have 
\[
\#  \Big\{\chi \pamod q\ \tmop{primitive}\,: \big|M_{\frac12}(\chi)\big| > \exp\big(\Lambda \sqrt{2 \log \log q}\big)\Big\} \ll qe^{-\Lambda ^2/36}.
\]

\end{lemma}
\begin{proof}
Applying Chebyshev's inequality and then the Cauchy-Schwarz inequality we have that 
\begin{align} \label{eq:moltailbd}
 & \#  \Big\{\chi \pamod q\ \tmop{primitive}:  \prod_{j=1}^J \big|M_{\frac12,j}(\chi)\big| \ge \exp\Big(\sqrt{\tfrac12 \log \log q}\Big) \Big\} \notag \\
 &  \qquad \quad
\le   e^{-( \frac12  \log \log q)^{1/2}} \sum_{\chi \pamod{q}} \prod_{j=1}^J \big|M_{\frac12,j}(\chi)\big|\nonumber \\
&  \qquad \quad \le   e^{-(\frac12  \log \log q)^{1/2}}\bigg( \sum_{\chi \pamod{q}} \prod_{j=1}^J \big|M_{\frac12,f,j}(\chi)\big|^2\bigg)^{1/2} \bigg( \sum_{\chi \pamod{q}} \prod_{j=1}^J \big|M_{\frac12,g,j}(\chi)\big|^2\bigg)^{1/2}\nonumber \\
 & \qquad \quad  \ll e^{-(\frac12  \log \log q)^{1/2}}  q (\log \log q)^{O(1)} \ll q e^{-(\frac13 \log \log q)^{1/2}},
\end{align}
where the second-to-last inequality follows from Lemma \ref{lem:secondmomentbd}.

Also, using \eqref{eq:ppower} we have that
\begin{equation} \label{eq:expandpower}
M_{\frac12,f,0}(\chi)=\sum_{\substack{p|n \Rightarrow p \in I_0 \\ \Omega(n) \le \ell_0}}\frac{  (-1/2)^{\Omega(n)}\nu(n)a_{f,J}(n) \chi(n)}{\sqrt{n}}=\sum_{\ell=0}^{\ell_0} \frac{(-\frac12 P_f(\chi))^{\ell}}{\ell!}.
\end{equation}
Additionally, for a complex number $z$ with $|z| \le N/e^2$ notice that
\begin{equation} \label{eq:exp}
\begin{split}
\sum_{n=0}^N \frac{z^n}{n!}&=e^z+O\Big( \sum_{n > N} \frac{(N/e^2)^n}{n!}\Big)=e^z\Big(1+O\big(e^{-N(1-1/e^2)}\big)\Big),
\end{split}
\end{equation}
where we used the inequality $n! >(n/e)^n$ in the last step. Hence, if $|P_f(\chi)|,|P_g(\chi)| \le \ell_0/e^2$, then it follows from \eqref{eq:expandpower} and \eqref{eq:exp} that
\begin{equation} \notag
M_{\frac12,0}(\chi)=e^{-\frac12  (P_f(\chi)+P_g(\chi))}\Big(1+O \big(e^{-\ell_0(1-1/e^2)}\big)\Big).
\end{equation}
Consequently,  if $|P_f(\chi)|,|P_g(\chi)| \le \ell_0/e^2$, then
\begin{equation} \label{eq:keybdM}
\big|M_{\frac12,0}(\chi)\big| \le e^{\frac12(|P_f(\chi)|+|P_g(\chi)|)}\big(1+o(1)\big).
\end{equation}

For each $f,g$ apply Lemma \ref{lem:devbd} with $V=\sqrt{\ell_0}$, so that $y^{V^2/9}=y^{\ell_0/9}<q$, to bound the number of $\chi \pmod q$ for which $|P_f(\chi)|+|P_g(\chi)|\ge 2V \sqrt{\frac12 \log \log y}$. Using this along with the estimate \eqref{eq:keybdM} applied to the characters $\chi \pmod q$ with 
\[
|P_f(\chi)|+|P_g(\chi)|\le 2V\sqrt{\tfrac12 \log \log y} \le \ell_0/e^2,
\] we have that 
\begin{align} \label{eq:zerobd}
& \#  \Big\{\chi \pamod q\ \tmop{primitive}:  \big|M_{\frac12,0}(\chi)\big| \ge \exp\Big(\Lambda \sqrt{\tfrac{1}{2} \log \log q}\Big)\Big\}\nonumber\\
&\qquad \ll  \#  \Big\{\chi \pamod q\ \tmop{primitive}:  |P_f(\chi)|+|P_g(\chi)| \ge \frac{\Lambda}{2} \sqrt{\tfrac{1}{2} \log \log q} \Big\} +qe^{-\ell_0/9}\nonumber \\
&\qquad \ll  q\big( e^{-\Lambda^2/36}+e^{-\ell_0/9}\big),
 \end{align}
using Lemma \ref{lem:devbd} for a second time, which is valid since $y^{\Lambda^2/36} < q$. 

To complete the proof note that if $|M_{\frac12,0}(\chi)|\ge \exp(\Lambda\sqrt{2 \log \log q})$ then \[\prod_{j=1}^J\big|M_{\frac12,j}(\chi)\big|\ge\exp\Big(\Lambda \sqrt{\tfrac{1}{2} \log \log q}\Big) \qquad \text{ or } \qquad \big|M_{\frac12,0}(\chi)\big|\ge\exp\Big(\Lambda \sqrt{\tfrac{1}{2} \log \log q}\Big).\] Hence, the result follows from \eqref{eq:moltailbd} and \eqref{eq:zerobd}. 
\end{proof}

\begin{proof}[Proof of Theorem \ref{thm2} assuming Propositions \ref{prop1}, \ref{prop2}, and \ref{prop3}]
We will only prove the first claim, as the second follows from a symmetric argument.
Let $I_{+}$ be as in \eqref{eq:Ipmdef}.
To shorten notation, write $L(\chi)=L(\tfrac12, f\otimes \chi)L(\tfrac12,g\otimes \overline \chi)$ and $W=\exp(\frac{c}{2} \sqrt{\log \log q})$. Observe for any $\delta>0$ we have that
\begin{equation} \label{eq:reduction}
\begin{split}
I_+=&\sum_{\substack{\chi \in S_+ \\ L(\chi)|M_{\frac12}(\chi)|^2>\delta}} L(\chi)\big|M_{\frac12}(\chi)\big|^2+O(q \delta)\\
=&\sum_{\substack{\chi \in S_+ \\ L(\chi)|M_{\frac12}(\chi)|^2>\delta \\ |M_{\frac12}(\chi)|^{2} \le W}} L(\chi)\big|M_{\frac12}(\chi)\big|^2+O\bigg( q\delta+\sum_{\substack{\chi \pamod q \\ |M_{\frac12}(\chi)|^{2}>W}} \big|L(\chi)M_{\frac12}(\chi)^2\big|\bigg).
\end{split}
\end{equation}
Applying the Cauchy-Schwarz inequality, we have by Proposition \ref{prop3} and Lemma \ref{lem:mollifierbd} that
\begin{align} \label{eq:boundobst}
 \bigg(\sum_{\substack{\chi \pamod q \\ |M_{\frac12}(\chi)|^{2}>W }} \big|L(\chi) M_{\frac12}(\chi)^2\big| \bigg)^2 & \le \#\Big\{\chi \pamod q : \big|M_{\frac12}(\chi)\big|^2 >W\Big\}  \sum_{\chi \pamod q } \big|L(\chi) M_{\frac12}(\chi)^2\big|^2 \nonumber\\
 & \ll q^2 e^{-c^2/288},
\end{align}
provided that $1 \le c \le 2\sqrt{2}(\log \log q)^{1/4}$.

 Recall that $I_+ \gg q$ by \eqref{estI+-}. We now take $\delta$ sufficiently small and fixed, and $c$ sufficiently large and fixed. Hence, also using \eqref{eq:reduction} and \eqref{eq:boundobst} it follows that
\[
 \sum_{\substack{\chi \in S_+ \\ L(\chi)>\delta W^{-1} }} L(\chi)\big|M_{\frac12}(\chi)\big|^2
 \ge \sum_{\substack{\chi \in S_+ \\ L(\chi)|M_{\frac12}(\chi)|^2>\delta \\ |M_{\frac12}(\chi)|^{2} \le W}} L(\chi)\big|M_{\frac12}(\chi)\big|^2
 \gg q.
\]
Therefore, applying the Cauchy-Schwarz inequality yields
\[
q \ll 
\# \Big\{ \chi \pamod q\ \tmop{primitive}: L(\chi)>\delta W^{-1} \Big
\}^{1/2}
\bigg(\ \,\sumstar_{\chi \pamod q}\big|L(\chi)M_{\frac12}(\chi)^2\big|^2 \bigg)^{1/2}. 
\]
Using Proposition \ref{prop3} and noting that $\delta W^{-1} \ge \exp(-c\sqrt{\log \log q})$ we conclude that
\[
\# \Big\{ \chi \pamod q\ \tmop{primitive}: L(\chi)>\exp\big(-c\sqrt{\log\log q}\big)\Big \} \gg q,
\]
as claimed.
\end{proof}

\section{Proofs of Proposition \ref{prop1} and Proposition \ref{prop2}} \label{sec:prop12}

We recall the twisted first moment of the product of twists of automorphic $L$-functions \cite[Theorem 5.1]{BFK}. To state the result, let $\varphi^{*}(q)$ denote the number of primitive characters modulo $q$, which is equal to $q-1$ for $q$ prime.

\begin{lemma}\label{BFKlemma} Let $n_1,n_2\in \mathbb N$ with $(n_1n_2/(n_1,n_2)^2,N)=1$. Then 
\[
\begin{split}
&\frac{1}{\varphi^{*}(q)}\ \sumstar_{\chi \pamod q} L(\tfrac12,f\otimes\chi)L(\tfrac12,g\otimes\overline{\chi})\chi(n_1) \overline{\chi(n_2)} =\frac12\sum_{m_1n_1=m_2n_2} \frac{\lambda_f(m_1)\lambda_g(m_2)}{\sqrt{m_1m_2}} V\Big( \frac{m_1m_2}{q^2N^2} \Big)\\
&\qquad+\frac{\varepsilon(f)\varepsilon(g)}{2}\sum_{m_1n_2=m_2n_1} \frac{\lambda_f(m_1)\lambda_g(m_2)}{\sqrt{m_1m_2}} V\Big( \frac{m_1m_2}{q^2N^2} \Big)+O_\varepsilon\left((n_1+n_2)^{3/2} q^{-1/144+\varepsilon}\right),
\end{split}
\]
where for $w>0$
\begin{equation}\label{Vformula}
V(w)= \frac{1}{2\pi i} \int_{(2)} \frac{L_{\infty}(s+\tfrac12,f)L_{\infty}(s+\tfrac12,g)}{L_{\infty}(\tfrac12,f)L_{\infty}(\tfrac12,g)} \cos\Big( \frac{\pi s}{12}\Big)^{-48} w^{-s} \, \frac{ds}{s} .
\end{equation}
\end{lemma}

\subsection{Main term estimates in terms of a random model}\label{subsectionrandom}
 Let $\{X(p)\}_p$ be i.i.d. random variables uniformly distributed on the unit circle. Additionally, let
\[
X(n)= \prod_{p^a || n} X(p)^a.
\]
 Following \cite[Sections 2, 3]{BELP}, we define the following random $L$-series, which streamline the proofs of Propositions \ref{prop1} and \ref{prop2}. First, let
\[
L(X)=\frac{1}{2} \Big( L^{f,g}(X)+\varepsilon(f)\varepsilon(g)\, L^{g,f}(X) \Big),
\]
where
\begin{align*}
L^{f,g}(X)&=\sum_{m_1,m_2 \ge 1} \frac{\lambda_f(m_1)\lambda_g(m_2)}{\sqrt{m_1m_2}} X(m_1) \overline{X(m_2)}\,V\Big(\frac{m_1m_2}{q^2N^2} \Big).
\end{align*}
Next, for $\tmop{Re}(s)>1$ define
\[
\begin{split}
L^{f,g}(s,X)&=\sum_{m,n \ge 1} \frac{\lambda_f(m)\lambda_g(n)}{(mn)^{s}} X(m) \overline{X(n)}=\prod_p \sum_{a,b \ge 0} \frac{\lambda_{f}(p^{a})\lambda_g(p^{b})}{p^{(a+b)s}} X(p)^{a} \overline{X(p)^{b}} \\
&= \prod_p  \prod_{j=1,2}\bigg(1-\frac{\alpha_{f,j}(p) X(p)}{p^{s}} \bigg)^{-1}\bigg(1-\frac{\alpha_{g,j}(p) \overline{X(p)}}{p^{s}} \bigg)^{-1}=:\prod_p L_p^{f,g}(s,X).
\end{split}
\]
Additionally, let
\[
L(s,X)=\frac{1}{2}\Big(L^{f,g}(s,X)+\varepsilon(f)\varepsilon(g) L^{g,f}(s,X)\Big).
\]

Finally, we introduce the random analogue of our mollifier. For each $j=0,\ldots,J$ define
\[
M_{\frac12,f,j}(X)=\sum_{\substack{p|n \Rightarrow p \in I_j \\ \Omega(n) \le \ell_j}}\frac{  (-1/2)^{\Omega(n)}\nu(n)a_{f,J}(n)  X(n)}{\sqrt{n}}
\]
where $J,I_j,\ell_j,\nu$, and $a_J$ are all as defined in Section \ref{sec:molandprop}. Additionally, let
\[
M_{\frac12,j}(X)=M_{\frac12,f,j}(X)M_{\frac12,g,j}(\overline{X})\qquad \text{and} \qquad M_{\frac12}(X)=\prod_{j=0}^JM_{\frac12,j}(X).
\]

 We shall prove the following propositions. 

\begin{proposition} \label{prop:random}
The functions $\mathbb E(L(s,X) M_{\frac12}(X)^2)$ and $\mathbb E(L(s,X) |M_{\frac12}(X)|^2)$ admit an analytic continuation to $\tmop{Re}(s)>0$. Furthermore, the following hold:
\begin{itemize}
\item[(i)]
\begin{align*}
\frac{1}{\varphi^{*}(q)}\ \sumstar_{\chi\pamod q} L(\tfrac12,f\otimes\chi)L(\tfrac12,g\otimes\overline{\chi})M_{\frac12,f}(\chi)^2M_{\frac12,g}(\overline{\chi})^2&=\mathbb E \Big( L(\tfrac12,X) M_{\frac12}(X)^2\Big)\\
&\qquad+O\left((\log q)^{-10}\right);
\end{align*}
\item[(ii)]
\begin{align*}
\frac{1}{\varphi^{*}(q)}\ \sumstar_{\chi\pamod q} L(\tfrac12,f\otimes\chi)L(\tfrac12,g\otimes\overline{\chi})\big|M_{\frac12,f}(\chi)M_{\frac12,g}(\overline{\chi})\big|^2&=\mathbb E \Big( L(\tfrac12,X) \big|M_{\frac12}(X)\big|^2\Big)\\
&\qquad+O\left((\log q)^{-10}\right).
\end{align*}
\end{itemize}
\end{proposition}

\begin{proposition} \label{prop:randombound}
The following hold:
\begin{itemize}

\item[(i)]
\[
\mathbb E \Big( L(\tfrac12,X)  M_{\frac12}(X)^2\Big) \asymp 1;
\]
\item[(ii)]
\[
\mathbb E \Big( L(\tfrac12,X)  \big|M_{\frac12}(X)\big|^2\Big) \ll (\log q)^{-1/2+\varepsilon}.
\]
\end{itemize}
\end{proposition}
\begin{proof}[Proof of Propositions \ref{prop1} and \ref{prop2} assuming Propositions \ref{prop:random} and \ref{prop:randombound}]
If $m,n \le q$, then by the orthogonality of characters we have that
\begin{equation} \label{eq:trivial}
\frac{1}{\varphi(q)}\sum_{\chi \pamod q} \chi(m) \overline {\chi(n)}=\mathbb E\Big(X(m) \overline {X(n)}\Big).
\end{equation}
Hence, from Lemma \ref{BFKlemma} and \eqref{eq:trivial}, we obtain  
\begin{align*}
\frac{1}{\varphi^{*}(q)}\ \sumstar_{\chi\pamod q} L(\tfrac12,f\otimes\chi)L(\tfrac12,g\otimes\overline{\chi})\chi(n_1) \overline{\chi(n_2)} &=\mathbb E \left( 
L(X) X(n_1) \overline{X(n_2)}
\right) \\
&\qquad+O_\varepsilon\left((n_1+n_2)^{3/2} q^{-1/144+\varepsilon} \right),
\end{align*}
where the contribution from $\chi= \chi_0$ is estimated as $O(1/q)$ since $q$ is prime. It follows that
\begin{align} \label{eq:twisted}
\frac{1}{\varphi^{*}(q)}\ \sumstar_{\chi\pamod q} L(\tfrac12,f\otimes\chi)L(\tfrac12,g\otimes\overline{\chi})M_{\frac12,f}(\chi)^2M_{\frac12,g}(\overline{\chi})^2& =\mathbb E \left( 
L(X)M_{\frac12}(X)^2
\right) \\
&\qquad+O(q^{-\delta})\nonumber,
\end{align}
for some $\delta>0$.
Applying Propositions \ref{prop:random}$(i)$ and \ref{prop:randombound}$(i)$ in \eqref{eq:twisted} completes the proof of Proposition \ref{prop1}.

Analogously, repeating the argument leading up to \eqref{eq:twisted} gives that
\begin{align} \label{eq:twisted2}
\frac{1}{\varphi^{*}(q)}\ \sumstar_{\chi\pamod q} L(\tfrac12,f\otimes\chi)L(\tfrac12,g\otimes\overline{\chi})\big|M_{\frac12,f}(\chi)M_{\frac12,g}(\overline{\chi})\big|^2& =\mathbb E \left( 
L(X)\big|M_{\frac12}(X)\big|^2
\right) \\
&\qquad+O(q^{-\delta}),\nonumber
\end{align}
for some $\delta>0$. Hence, using Propositions \ref{prop:random}$(ii)$ and \ref{prop:randombound}$(ii)$ in \eqref{eq:twisted2} completes the proof of Proposition \ref{prop2}.
\end{proof}

\subsection{Preparatory lemmas}\label{sec:prime}
The arguments of this section follow the approach of \cite[Section 3]{BELP} and we will reproduce the details here with appropriate modifications to make the presentation more self contained. Since the arguments for each of Proposition \ref{prop:random}$(i)$,$(ii)$ and Proposition \ref{prop:randombound}$(i)$,$(ii)$ are similar we will focus on the results needed for Propositions \ref{prop:random}$(ii)$, \ref{prop:randombound}$(ii)$ since these are further from \cite{BELP}. The proof of Proposition \ref{prop:random} uses a contour shifting argument, which is justified from the analytic continuation of and bounds for the expected values of each of $L(s,X)M_{\frac12}(X)^2$ and $L(s,X)|M_{\frac12}(X)|^2$, which we establish in this section. 

Recall $y=q^{\theta_0}$ and $x=q^{\theta_J}$, as defined in Section \ref{sec:molandprop}. To state the next result, let
\[
\mathfrak C=N^2(|s+2|+1)^2(|s|+1)^2
\]
denote the analytic conductor of $f \otimes g$.
By well-known properties of Rankin-Selberg $L$-functions, see \cite[Lemma 3.6]{BELP} for a proof, we have there exists $c>0$ such that for  $-c/\log \mathfrak C \le \tmop{Re}(s) \le 2$, $|\tmop{Im}(s)| \le e^{\sqrt{\log q}}$, and $1 \le j \le J$ that
    \begin{equation} \label{eq:RSbds}
     \sum_{p \in I_j} \frac{\lambda_f(p)\lambda_g(p)}{p^{2s+1}}=O(1), \qquad \sum_{p \in I_j} \frac{ w_J(p)\lambda_f(p)\lambda_g(p)}{p^{s+1}} =O(1).
     \end{equation}

\begin{lemma} \label{lem:mediumprimeest}
     Let $\sigma_0=\max\{ -\tmop{Re}(s),\frac{1}{\log q}\}$.
There exists $c>0$ such that for
$ -\frac{c}{\log \mathfrak C} \le \tmop{Re}(s) \le 2$ and $|\tmop{Im}(s)| \le e^{\sqrt{\log q}}$ the following hold:
\begin{itemize}
\item[(i)]\begin{equation} \notag
\mathbb E\bigg( \prod_{y< p \le x} L_p^{f,g}(s+\tfrac12,X) \prod_{j=1}^J M_{\frac12,j}(X)^2 \bigg)  \ll (\log \log q)^{O(1)};
\end{equation}
\item[(ii)]
   \begin{equation} \notag
\mathbb E\bigg( \prod_{y< p \le x} L_p^{g,f}(s+\tfrac12,X) \prod_{j=1}^J M_{\frac12,j}(X) ^2\bigg)  \ll (\log \log q)^{O(1)} q^{41 \eta^{1/4} \sigma_0};
\end{equation}
\item[(iii)]   \begin{equation} \notag
\mathbb E\bigg( \prod_{y< p \le x} L_p^{f,g}(s+\tfrac12,X) \prod_{j=1}^J \big|M_{\frac12,j}(X)\big| ^2\bigg)  \ll (\log \log q)^{O(1)} q^{41 \eta^{1/4} \sigma_0}.
\end{equation}
\end{itemize}
\end{lemma}
    \begin{proof}
    The arguments used to prove $(i)-(iii)$ are similar and we will provide the full details of $(iii)$ and omit the arguments used to prove $(i)$ and $(ii)$.
    Let $\beta_f(n)=\lambda_f(n)n^{-s}$ and $\gamma_f(n)=(-1/2)^{\Omega(n)}\nu(n)a_{f,J}(n)$. Clearly, $\beta_f,\gamma_f$ are multiplicative functions. Then for each $1 \le j \le J$ we have
    \begin{align} \label{eq:expectedvalue}
  & \mathbb E\bigg(   \big|M_{\frac12,j}(X)\big| ^2\prod_{p \in I_j} L_p^{f,g}(s+\tfrac12,X) \bigg) \nonumber\\
    &\qquad=\mathbb E\bigg(  \sum_{p|m_1m_2 \Rightarrow p\in I_j} \frac{\beta_f(m_1)\beta_g(m_2) X(m_1) \overline {X(m_2)}}{ \sqrt{m_1m_2}}\nonumber\\
&\qquad\qquad\qquad\times\sum_{\substack{p|n_1n_2n_3n_4 \Rightarrow p \in I_j \\ \Omega(n_i)\le \ell_j,1\le i\le 4}} \frac{\gamma_f(n_1)\gamma_g(n_2)\gamma_f(n_3)\gamma_g(n_4)X(n_1n_4) \overline{X(n_2n_3)}}{\sqrt{n_1n_2n_3n_4}}\bigg)\nonumber \\
    &\qquad=\sum_{\substack{m_1,m_2, n_1,n_2,n_3,n_4 \\ p|m_1n_1n_4 \Rightarrow p \in I_j \\ m_1n_1n_4=m_2n_2n_3 \\ \Omega(n_i) \le \ell_j,1\le i\le 4}} \frac{ \beta_f(m_1) \beta_g(m_2)\gamma_f(n_1)\gamma_g(n_2)\gamma_f(n_3)\gamma_g(n_4)}{m_1n_1n_4}.
    \end{align}
   For any $n, \ell \in \mathbb N$ we have
    \[
    \mathbf{1}_{\Omega(n)=\ell}= \frac{1}{2\pi i} \int_{|z|=2} z^{\Omega(n)-\ell} \frac{dz}{z},
    \]
    and hence
    \[
    \mathbf{1}_{\Omega(n) \le \ell_j} =\frac{1}{2\pi i} \int_{|z|=2} z^{\Omega(n)} \frac{1-z^{-\ell_j-1}}{1-z^{-1}} \frac{dz}{z}.
    \]
    Therefore,  the left-hand side of \eqref{eq:expectedvalue} equals
    \begin{equation}\label{eq:expanded}
    \frac{1}{(2\pi i)^4} \int_{\substack{|z_i|=2\\1\le i\le 4}}\prod_{i=1}^{4}\frac{1-z_i^{-\ell_j-1}}{1-z_i^{-1}} \Sigma(\textbf{z}) \frac{dz_1dz_2dz_3dz_4}{z_1z_2z_3z_4},
    \end{equation}
    where
    \[
    \Sigma(\textbf{z})=\sum_{\substack{m_1,m_2, n_1,n_2,n_3,n_4 \\ p|m_1n_1n_4 \Rightarrow p \in I_j \\ m_1n_1n_4=m_2n_2n_3 }} \frac{ \beta_f(m_1) \beta_g(m_2)\gamma_f(n_1)\gamma_g(n_2)\gamma_f(n_3)\gamma_g(n_4)\prod_{i=1}^4z_i^{\Omega(n_i)}}{m_1n_1n_4},
    \]
   which can be seen to be absolutely convergent from the analysis below.
      
    We write $\gamma_{z,f}(n)=z^{\Omega(n)}\gamma_f(n)$ and $n=m_1n_1n_4$ to see that 
    \begin{equation} \label{eq:eulerexpand}
    \begin{split}
     \Sigma(\textbf{z})=& \sum_{p|n \Rightarrow p \in I_j} \frac{(\beta_f\ast \gamma_{z_1,f}\ast \gamma_{z_4,g})(n) (\beta_g\ast \gamma_{z_2,g}\ast \gamma_{z_3,f})(n)}{n} \\
     =& \prod_{p \in I_j} \sum_{a\geq 0} \frac{(\beta_f\ast \gamma_{z_1,f}\ast \gamma_{z_4,g})(p^a)(\beta_g\ast \gamma_{z_2,g}\ast \gamma_{z_3,f})(p^a)}{p^a} \\
     =& \prod_{p \in I_j} \bigg(1+\frac{1}{p}\Big(\lambda_f(p)p^{-s}-w_J(p)\big(\lambda_f(p)z_1+\lambda_g(p)z_4\big)/2\Big) \\
&\qquad\qquad\times \Big(\lambda_g(p)p^{-s}-w_J(p)\big(\lambda_g(p)z_2+\lambda_f(p)z_3\big)/2 \Big)+O\Big( \frac{1+p^{-4\text{Re}(s)}}{p^2}\Big) \bigg).
     \end{split}
     \end{equation}
 Using \eqref{eq:RSbds} and that \[\sum_{p \in I_j} \frac{\lambda_f(p)^2w_J(p)^2}{p} \ll 1,\]
which follows from Deligne's bound, we have that
    \begin{align*}
    \Sigma(\textbf{z})=\exp\bigg(-\frac{z_3}{2}\sum_{p \in I_j} \frac{\lambda_f(p)^2 w_J(p)}{p^{1+s}}-\frac{z_4}{2}\sum_{p \in I_j} \frac{\lambda_g(p)^2 w_J(p)}{p^{1+s}}\bigg)H(\textbf{z}),
    \end{align*}
    where $H(\textbf{z})$ is an analytic function for $|z_i| \le 2$, $1\le i\le4$, with $|H(\textbf{z})| \ll 1$ uniformly for $s$ satisfying the hypotheses of the lemma.

     We now split the proof into two cases. First consider the case, $\theta_j < \frac{1}{\sigma_0 \log q}$. Then we have by Deligne's bound that     \begin{equation} \label{eq:smallj}
     \sum_{p \in I_j} \frac{\lambda_f(p)^2 w_J(p)}{p^{1+s}} \ll e^{\sigma_0 \theta_j \log q} \sum_{p \in I_j} \frac{1}{p} \ll 1,
     \end{equation}
     so that $\Sigma(\textbf{z}) \ll1$ in this case and the left-hand side of \eqref{eq:expectedvalue} is $\ll 1$.
    
    It remains to consider the case that $\theta_j \ge \frac{1}{\sigma_0 \log q}$. Using that $\Sigma(\textbf{z})$ is analytic in each variable, expanding $(1-z_i^{-\ell_j-1})/(1-z_i)$, $1\le i\le 4$, as geometric sums, and using Cauchy's integral formula we get that the left-hand side of \eqref{eq:expectedvalue} equals
    \begin{equation} \label{eq:derivrep}
    \sum_{\substack{0 \le k_i \le \ell_j\\1\le i\le4}} \frac{1}{k_1!k_2!k_3!k_4!} \frac{\partial^{k_1+k_2+k_3+k_4}}{\partial z_1^{k_1} \partial z_2^{k_2}\partial z_3^{k_3} \partial z_4^{k_4}}  \Sigma(\textbf{z}) \bigg|_{\textbf{z}=\textbf{0}}.
    \end{equation}
   Using Cauchy's integral formula once again, we see that $$\frac{\partial^{k_1+k_2+k_3+k_4}}{\partial z_1^{k_1} \partial z_2^{k_2}\partial z_3^{k_3} \partial z_4^{k_4}} H(\textbf{z}) \bigg|_{\textbf{z}=\textbf{0}}\ll \frac{k_1! k_2! k_3!k_4!}{2^{k_1+k_2+k_3+k_4}}.$$
Thus we can bound the expression in \eqref{eq:derivrep} by 
    \begin{equation}\label{bddev}
    \ll 4^{4\ell_j} \bigg( \sum_{p \in I_j} \frac{\lambda_g(p)^2w_J(p) }{p^{1-\sigma_0}} \bigg)^{2\ell_j} \bigg( \sum_{p \in I_j} \frac{\lambda_f(p)^2w_J(p)}{p^{1-\sigma_0}} \bigg)^{2\ell_j},
    \end{equation}
    where the term $4^{4\ell_j}$ comes from applying the product rule.
    Using Deligne's bound and that $\ell_j \le 2\theta_j^{-3/4}$ we conclude that the left-hand side of \eqref{eq:expectedvalue} is 
    \begin{equation} \label{eq:largej}
    \ll \bigg(16 \sum_{p \in I_j} \frac{1}{p^{1-\sigma_0}}\bigg)^{4\ell_j} \le 18^{4\ell_j} e^{8 \theta_j^{1/4} \sigma_0 \log q}.
    \end{equation}
    
    Write $J_1$ for the largest $j=1,\ldots,J$ such that $\theta_j <\frac{1}{\sigma_0 \log q} $. Using \eqref{eq:smallj}, \eqref{eq:largej} and that for $J_1 \le j \le J$, $\ell_j \le 2 (\sigma_0 \log q)^{3/4}$ we have that the left hand side of $(iii)$ is
    \[
   \ll (\log \log q)^{O(1)} \prod_{J_1 \le j \le J}  18^{4\ell_j} e^{8 \theta_j^{1/4} \sigma_0 \log q} \ll (\log \log q)^{O(1)} e^{O((\sigma_0 \log q)^{3/4})}\,  e^{40 \eta^{1/4} \sigma_0 \log q}.
    \]
    This finishes the proof.
     \end{proof}

The pieces of our mollifiers corresponding to $j=0$ behave very much like Euler products, which allow for direct computation. We quantify this relation in the following lemma, and to state it let 
     \[
    \widetilde M_{\frac12,0}(X)=\sum_{\substack{p|n \Rightarrow p \in I_0}} \frac{(-1/2)^{\Omega(n)}}{\sqrt{n}} (\nu a_{f,J} X \ast \nu a_{g,J}\overline X)(n)= \prod_{p \in I_0} \widetilde M_{\frac12,p}(X),
    \]
    where
    \begin{equation} \label{eq:Mtildedef}
    \widetilde M_{\frac12,p}(X)=\widetilde M_{\frac12,p,f}(X) \widetilde M_{\frac12,p,g}(\overline X), \qquad \widetilde M_{\frac12,p,f}(X)=\sum_{k\geq0} \frac{(-1/2)^ka_{f,J}(p)^k}{k!p^{k/2}} X(p)^k.
    \end{equation}
    
\begin{lemma} \label{lem:multiplicative} 
    Suppose  $\tmop{Re}(s) \ge -\frac{( \log \log q)^2}{\log q}$. The following hold.
    \begin{itemize} 
    \item[(i)] For $(f_1,f_2)=(f,g)$ or $(f_1,f_2)=(g,f)$ we have
    \[
    \begin{split}
    &\mathbb E\bigg( M_{\frac12,0}(X)^2\prod_{p\in I_0}L_p^{f_1,f_2}(s+\tfrac12,X)  \bigg)= \mathbb E\bigg(\widetilde M_{\frac12,0}(X)^2 \prod_{p\in I_0}L_p^{f_1,
    f_2}(s+\tfrac12,X)  \bigg)+O\left( (\log q)^{-10}\right).
    \end{split}
    \]   
    \item[(ii)] We have
    \[
    \begin{split}
    &\mathbb E\bigg( \big|M_{\frac12,0}(X)\big|^2\prod_{p\in I_0}L_p^{f,g}(s+\tfrac12,X)  \bigg)= \mathbb E\bigg(\big|\widetilde M_{\frac12,0}(X)\big|^2 \prod_{p\in I_0}L_p^{f,
    g}(s+\tfrac12,X)  \bigg)+O\left( (\log q)^{-10}\right).
    \end{split}
    \]   
    \end{itemize}
    \end{lemma}
    \begin{proof}
    We will only prove $(ii)$; the proof of $(i)$ follows from a similar argument.
    Recall that  $\gamma_f(n)=(-1/2)^{\Omega(n)}\nu(n)a_{f,J}(n)$ and define
    \begin{align*}
    R(X)&:=\big|\widetilde M_{\frac12,0}(X)\big|^2-\big|M_{\frac12,0}(X)\big|^2\\
&=\sum_{\substack{p|m_1m_2n_1n_2 \Rightarrow p \in I_0 \\ \max\{ \Omega(m_1),\Omega(m_2),\Omega(n_1),\Omega(n_2)\} > \ell_0}}  \frac{\gamma_f(m_1)\gamma_f(m_2)\gamma_g(n_1)\gamma_g(n_2)}{\sqrt{m_1m_2n_1n_2}} X(m_1n_2)\overline {X(m_2n_1)}.
    \end{align*}
    Writing $L_0^{f,g}(s,X)=\prod_{p \in I_0} L_p^{f,g}(s,X)$ and using the Cauchy-Schwarz inequality we have
    \begin{align*}
\mathbb E\Big( \big|L_0^{f,g}(s+\tfrac12,X) R(X) \big|\Big)^2
  & \le    \mathbb E \Big(|L_0^{f,
   g}(s+\tfrac12,X) |^2 \Big) \mathbb E\big(|R(X)|^2\big)\\
&\le (\log q)^{O(1)} \mathbb E\big(|R(X)|^2\big).
    \end{align*}
Here the last inequality uses the bound
\[
\mathbb E \Big(|L_0^{f,
   g}(s+\tfrac12,X) |^2 \Big) =\prod_{p \in I_0}\bigg( 1+O\Big(\frac{1}{p^{1+2\tmop{Re}(s)}}\Big)\bigg) =\prod_{p \in I_0}\bigg( 1+O\Big(\frac{1}{p}\Big)\bigg) \ll (\log q)^{O(1)}
\]
obtained in \cite[Eq'ns (3.12)-(3.13)]{BELP}, where we have used that $p^{-\tmop{Re}(s)}\ll 1$ for $\tmop{Re}(s) > \frac{-(\log \log q)^2}{\log q}$ and $p \in I_0$.
 
        We next estimate $\mathbb E( |R(X)|^2)$, which equals
    \begin{align*}
    &\sum_{\substack{p|m_1m_2m_3m_4n_1n_2n_3n_4 \Rightarrow p \in I_0 \\ \max\{ \Omega(m_1),\Omega(m_2),\Omega(n_1),\Omega(n_2)\} > \ell_0 \\ \max\{ \Omega(m_3),\Omega(m_4),\Omega(n_3),\Omega(n_4)\} > \ell_0}} \frac{\prod_{j=1}^4\gamma_f(m_j)\gamma_g(n_j)}{\sqrt{\prod_{j=1}^4m_jn_j}}\mathbb E\Big(X(m_1m_4n_2n_3)\overline {X(m_2m_3n_1n_4)}\Big).
    \end{align*}
 We write $n=m_1m_4n_2n_3=m_2m_3n_1n_4$. Using the fact that $$\max\{ \Omega(m_1),\Omega(m_2),\Omega(n_1),\Omega(n_2)\} > \ell_0\quad \text{and}\quad \max\{ \Omega(m_3),\Omega(m_4),\Omega(n_3),\Omega(n_4)\} > \ell_0$$ imply $2^{\Omega(n^4)-\ell_0} \ge 1$, we get upon applying Deligne's bound  $|a_{f,J}(n)|  \le 2^{\Omega(n)}$ that  
    \begin{equation}\label{eq:Rerrorbd}
    \mathbb E\big(|R(X)|^2\big) \le \frac{1}{2^{\ell_0}} \sum_{p|n \Rightarrow p \in I_0 } \frac{C^{\Omega(n)}}{n}   \ll \frac{(\log q)^{O(1)}}{2^{\ell_0}},
    \end{equation} for some $C>0$. This completes the proof.
     \end{proof}
     
     In the next lemma, we explicitly evaluate the Euler product factors that arise from Lemma \ref{lem:multiplicative}.
    To state the result, recall the definition of $\widetilde M_{\frac12,p}(X)$ given in \eqref{eq:Mtildedef} and let
    \[
    \mathcal G_p^{f,g}(s)=\mathbb E \left(L_p^{f,g}(s+\tfrac12, X) \widetilde M_{\frac12,p}(X) ^2\right)
    \]
    and
    \[
    \mathcal H_p^{f,g}(s)=\mathbb E \left(L_p^{f,g}(s+\tfrac12, X) \big|\widetilde M_{\frac12,p}(X)\big| ^2\right).
    \]
    
    \begin{lemma} \label{lem:gfest} For $p^{-\tmop{Re}(s)} \le 2$ the following hold:
    \begin{itemize}
    \item[(i)] \[
        \mathcal G_p^{f,g}(s)=1+\frac{(p^{-s}-w_J(p))^2\lambda_f(p)\lambda_g(p)}{p} +O\Big( \frac{1}{p^2}\Big);
        \]
        \item[(ii)]
        \[\mathcal G_p^{g,f}(s)=1+\frac{(\lambda_f(p)p^{-s}-w_J(p)\lambda_g(p))(\lambda_g(p)p^{-s}-w_J(p)\lambda_f(p))}{p} +O\Big( \frac{1}{p^2}\Big);
        \]
        \item[(iii)]
      \begin{align*} \mathcal H_p^{f,g}(s)&=1+\frac{\big(\lambda_f(p)p^{-s}-w_J(p)(\lambda_f(p)+\lambda_g(p))/2\big)\big(\lambda_g(p)p^{-s}-w_J(p)(\lambda_f(p)+\lambda_g(p))/2\big)}{p}\\
      &\qquad+O\Big( \frac{1}{p^2}\Big).
        \end{align*}
        \end{itemize}
    \end{lemma}
    \begin{proof}
    The proofs of $(i)-(iii)$ are similar and we will only prove $(iii)$. Observe that 
      for $\tmop{Re}(s)>0$ and $k$
a positive integer 
\[
\bigg(\sum_{k\ge 0} \frac{\lambda_{f}(p^{k} ) X(p)^{k}}{p^{k s}} \bigg) \,  \widetilde M_{\frac12,p,f}(X) \widetilde M_{\frac12,p,g}(X) =\sum_{k \ge 0} n_p^{f,g}(s,k) X(p)^k.
\]
where $\widetilde M_{\frac12,p,f}(X)$ is as given in \eqref{eq:Mtildedef} and
\[
n_p^{f,g}(s,k)=\sum_{\substack{k_1+k_2+k_3=k \\ k_1, k_2 ,k_3 \ge 0}} \frac{\lambda_{f}(p^{k_1})(-1/2)^{k_2+k_3} a_{f,J}(p)^{k_2}a_{g,J}(p)^{k_3}} {p^{k_1 s+(k_2+k_3)/2}k_2!k_3!}.
\]
Hence
    \begin{align} \label{eq:expectation1}
    \mathbb E\left(L_p^{f,g}(s+\tfrac12, X) |\widetilde M_{\frac12,p}(X)|^2 \right)=&\sum_{k_1,k_2 \ge 0} n_p^{f,g}(s+\tfrac12,k_1) n_p^{g,f}(s+\tfrac12,k_2) \mathbb E\left( X(p)^{k_1} \overline {X(p)^{k_2}} \right)  \nonumber \\
    =& \sum_{k\geq 0} n_p^{f,g}(s+\tfrac12,k)n_p^{g,f}(s+\tfrac12,k).
    \end{align}

    We have that $n_p^{f,g}(s+\tfrac12,0)=1$ and 
    \begin{equation} \label{eq:npfg1}
    n_p^{f,g}(s+\tfrac12,1)=\frac{(\lambda_{f}(p)p^{-s}-w_J(p) (\lambda_{f}(p)+\lambda_g(p))/2)}{\sqrt{p}}.
    \end{equation}
    Additionally, for each $k \ge 1$ using Deligne's bound we get that
    \begin{equation} \label{eq:npfg2}
    n_p^{f,g}(s+\tfrac12,k) \ll \frac{2^k}{p^{k/2}} \max\{1,p^{-\tmop{Re}(s)k}\}.
    \end{equation}
    Applying \eqref{eq:npfg1} and \eqref{eq:npfg2} in \eqref{eq:expectation1} completes the proof.
    \end{proof}

    Finally, we require Lemma 3.7 from \cite{BELP}, which we now state.
    \begin{lemma} \label{lem:largeprimeest}
There exists $c>0$ such that 
for $-\frac{c}{\log \mathfrak C} \le \tmop{Re}(s) \le 2$, $z  \ge \exp(\log \mathfrak C \log \log \mathfrak C)$, and  for $(f_1,f_2)=(f,g)$ or $(f_1,f_2)=(g,f)$ we have 
\[
\mathbb E\bigg(\prod_{p > z} L_p^{f_1,f_2}(s+\tfrac12,X) \bigg)=1+O\left(z^{-c/\log \mathfrak C}(\log z)^{2}\right).
\]
\end{lemma}

 \subsection{Proof of Propositions \ref{prop:random} and \ref{prop:randombound}}
 Let us record a few useful consequences from Section \ref{sec:prime}. By Lemmas \ref{lem:multiplicative} and \ref{lem:gfest}  we have for $(f_1,f_2)=(f,g)$ or $(f_1,f_2)=(g,f)$ that 
  \begin{equation} \label{eq:useful}
  \begin{split}
&\mathbb E\bigg( M_{\frac12,0}(X)^2 \prod_{p\in I_0}L_p^{f_1,f_2}(s+\tfrac12,X)  \bigg) \ll (\log q)^{O(1)},\\
&\mathbb E\bigg( \big|M_{\frac12,0}(X)\big|^2 \prod_{p\in I_0}L_p^{f,g}(s+\tfrac12,X)  \bigg) \ll (\log q)^{O(1)}
\end{split}
\end{equation}
  for $\tmop{Re}(s) \ge -\frac{(\log \log q)^2}{\log q}$. 
  Hence for $ -\frac{(\log \log q)^2}{\log q} \le \tmop{Re}(s)  \le 2$ and $|\tmop{Im}(s)| \le e^{\sqrt{\log q}}$, Lemmas \ref{lem:mediumprimeest} and \ref{lem:largeprimeest} along with \eqref{eq:useful} give
  \begin{equation} \label{eq:fupperbd1}
  \mathbb E\Big( L(s+\tfrac12,X) M_{\frac12}(X)^2  \Big) \ll (\log q)^{O(1)}+ (\log q)^{O(1)}e^{41 \eta^{1/4} (\log \log q)^2}
  \end{equation}
  and
   \begin{equation} \label{eq:fupperbd12}
  \mathbb E\Big( L(s+\tfrac12,X) \big|M_{\frac12}(X)\big|^2  \Big) \ll  (\log q)^{O(1)}e^{41 \eta^{1/4} (\log \log q)^2}.
  \end{equation}

 \begin{proof}[Proof of Proposition \ref{prop:random}]
 For $\tmop{Re}(s)>\tfrac12 $ define
  \begin{equation} \notag
   F_1(s)=\mathbb E \Big(L(s+\tfrac12,X)M_{\frac12}(X)^2 \Big)
   \end{equation}
   and
 \begin{equation} \notag
   F_2(s)=\mathbb E \Big(L(s+\tfrac12,X)\big|M_{\frac12}(X)\big|^2 \Big).
   \end{equation}
   For $(f_1,f_2)=(f,g)$ and $(f_1,f_2)=(g,f)$ \cite[Eq'n (3.2)]{BELP}, establishes the analytic continuation of
   \[
   \mathbb E \bigg( \prod_{p>x} L_p^{f_1,f_2}(s+\tfrac12,X)  \bigg)
   \]
   to $\tmop{Re}(s)>-\tfrac12$. Hence, each $F_1(s)$ and $F_2(s)$
   admits an analytic continuation to $\tmop{Re}(s)>-\tfrac12$. 
   
   Moreover, by \eqref{eq:fupperbd1} and \eqref{eq:fupperbd12} we have that 
     \begin{equation} \label{eq:fbd1}
   |F_1(s)|,|F_2(s)|  \ll e^{42 \eta^{1/4} (\log \log q)^2}
   \end{equation}
   uniformly for $ - \frac{(\log \log q)^2}{ \log q} \le \tmop{Re}(s) \le 2$ and $|\tmop{Im}(s)| \le B \log q$ for any fixed $B>0$. Additionally, we have  
   \begin{equation}\label{eq:fbd2}
   |F_1(s)|,|F_2(s)|  \ll (\log q)^{O(1)}\
   \end{equation}
   in the region $\tmop{Re}(s) \ge \tfrac12+\varepsilon$ since
   \[
   \mathbb E \Big( \big|M_{\frac12}(X)\big|^2\Big) \ll (\log q)^{O(1)},
   \]
which follows from a straightforward modification of the proof of Lemma \ref{lem:secondmomentbd} that we will omit.
   
   Applying the Mellin inversion theorem we see that
   \begin{equation} \label{eq:mellin}
   \begin{split}
   &\mathbb E \Big(L(X)M_{\frac12}(X)^2 \Big) = \frac{1}{2\pi i} \int_{(2)}   \frac{L_{\infty}(s+\tfrac12,f)L_{\infty}(s+\tfrac12,g)}{L_{\infty}(\tfrac12,f)L_{\infty}(\tfrac12,g)}\cos\Big( \frac{\pi s}{12}\Big)^{-48}  \left( q^2 N^2 \right)^sF_1(s)  \, \frac{ds}{s}.
   \end{split}
   \end{equation}
   Since for fixed $\sigma>0$, Stirling's formula gives that $|\Gamma(\sigma+it)|\ll (|t|+1)^{\sigma-\frac12} e^{-\pi|t|}$, by \eqref{eq:fbd2} we may truncate the integral in \eqref{eq:mellin} to $|\tmop{Im}(s) |\le B \log q$ at the cost of an error term of size $O(q^{-1})$ where $B$ is a sufficiently large absolute constant. We now shift contours to $\tmop{Re}(s)=-\frac{(\log \log q)^2}{\log q}$,  pick up a simple pole at $s=0$, estimate the horizontal and left contours using the bound \eqref{eq:fbd1} to get that
   \[
   \begin{split}
  & \mathbb E \Big(L(X)M_{\frac12}(X)^2 \Big)=F_1(0)+O\left((\log q)^{-10}\right).
   \end{split}
   \]
   which completes the proof of $(i)$.

   Similarly we have
   \[
   \begin{split}
  & \mathbb E \Big(L(X)\big|M_{\frac12}(X)\big|^2 \Big)=F_2(0)+O\left((\log q)^{-10}\right).
   \end{split}
   \]
   and this completes the proof of $(ii)$.
\end{proof}

    \begin{proof}[Proof of Proposition \ref{prop:randombound}]
    We will first prove $(i)$.
        Using \eqref{eq:expectedvalue} with $s=0$ and recalling $\gamma_f(n)=(-1/2)^{\Omega(n)} \nu(n)a_{f,J}(n)$,
    we have for each $0\leq j\leq J$ that
    \[
       \mathbb E \bigg( M_{\frac12,j}(X)^2 \prod_{p \in I_j} L_p^{f,g}(\tfrac12, X) \bigg)= \sum_{\substack{m_1,m_2,n_1,n_2 ,n_3,n_4\\ p |m_1n_1n_3 \Rightarrow p \in I_j \\ m_1n_1n_3=m_2n_2n_4 \\ \Omega(n_i) \le \ell_j,1\le i\le 4}} \frac{\lambda_f(m_1)\lambda_g(m_2)\gamma_f(n_1)\gamma_g(n_2)\gamma_f(n_3)\gamma_g(n_4) }{m_1 n_1n_3}.
       \]
       We now wish to remove the condition $\Omega(n_i) \le \ell_j$ so that we can express the sum as an Euler product. Arguing as in \eqref{eq:Rerrorbd} we see that there exists $C>0$ such that the sum on the right-hand side is
       \begin{equation*}
       \begin{split}
       =& \sum_{\substack{ p |n \Rightarrow p \in I_j  }} \frac{(\lambda_f\ast \gamma_f\ast \gamma_f)(n)(\lambda_g \ast \gamma_g\ast \gamma_g)(n) }{m} +O\bigg(\frac{1}{2^{\ell_j}} \sum_{p | n \Rightarrow p \in I_j} \frac{C^{\Omega(n)}}{n}\bigg) \\
       =&\prod_{p \in I_j}\left(1+\frac{(1-w_J(p))^2 \lambda_f(p) \lambda_g(p)}{p}+O\Big(\frac{1}{p^2} \Big) \right)+O \left( \frac{\mathbf{1}_{j=0} (\log q)^{O(1)}+1}{2^{\ell_j}}\right),
       \end{split}
       \end{equation*}
       which is the same as (3.23) of \cite{BELP}.
The arguments in Section 3.5 of \cite{BELP} then yield 
\begin{equation}\label{1stest}
\mathbb E\Big( L^{f,g}(\tfrac12,X)M_{\frac12}(X)^2\Big) \asymp 1.
\end{equation}

For $L^{g,f}(\tfrac12,X)M_{\frac12}(X)^2$, we use Lemma \ref{lem:largeprimeest} to estimate the contribution of the primes $p>x$ and Lemmas \ref{lem:multiplicative}$(i)$ and \ref{lem:gfest}$(ii)$ to estimate that of the primes $c_0 < p \le y$ to get 
        \begin{align}\label{2ndpart}
  & \mathbb E\Big( L^{g,f}(\tfrac12,X)M_{\frac12}(X)^2\Big)\nonumber\\
  &\ = \big(1+o(1)\big)  \prod_{p \le c_0}\mathbb E\Big( L_p^{g,f}(\tfrac12,X)\Big) \prod_{y< p \le x} \mathbb E\bigg(  L_p^{g,f}(\tfrac12,X) \prod_{j=1}^J M_{\frac12,j}(X)^2\Big) \bigg) \\
   &\quad\times   \bigg( 
\prod_{c_0 < p \le y}  \left(1+\frac{(\lambda_f(p)-w_J(p)\lambda_g(p)) (\lambda_g(p)-w_J(p)\lambda_f(p))}{p}+O\Big(\frac{1}{p^2} \Big)\right) +O\big((\log q)^{-10}\big)\bigg).\nonumber
    \end{align}
    The product over $p \le c_0$ is $\ll 1$ by \cite[Lemma 3.4]{BELP}. By Lemma \ref{lem:mediumprimeest}$(ii)$, the product over $y< p \le x$ is $\ll (\log \log q)^{O(1)}$.
   For the product over $c_0< p \le y$, we use the Prime Number Theorem for Rankin-Selberg $L$-functions, which implies \begin{equation}\label{RSexp}
   \sum_{c_0 < p \le y} \frac{\lambda_f(p)\lambda_g(p)}{p}=O(1) \qquad  \text{ and } \qquad \sum_{c_0 < p \le y} \frac{\lambda_f(p)^2}{p}=\log \log y+O(1),
   \end{equation}
   to see that it is $\ll (\log y)^{-2}$.
   Thus,
\[
\mathbb E\Big( L^{g,f}(\tfrac12,X)M_{\frac12}(X)^2\Big)=O_\varepsilon\big( (\log q)^{-2+\varepsilon}\big).
\]
Combining the preceding estimate with \eqref{1stest} we obtain $(i)$.

For $(ii)$, similar to \eqref{2ndpart} we have
        \begin{align*}
  & \mathbb E\Big( L^{f,g}(\tfrac12,X)\big|M_{\frac12}(X)\big|^2\Big)\nonumber\\
  &\ = \big(1+o(1)\big)  \prod_{p \le c_0}\mathbb E\Big( L_p^{f,g}(\tfrac12,X)\Big) \prod_{y< p \le x} \mathbb E\bigg(  L_p^{f,g}(\tfrac12,X) \prod_{j=1}^J \big|M_{\frac12,j}(X)\big|^2 \bigg) \\
   &\quad\times   \bigg( 
\prod_{c_0 < p \le y}  \bigg(1+\frac{\big(\lambda_f(p)-w_J(p)(\lambda_f(p)+\lambda_g(p))/2\big)\big(\lambda_g(p)-w_J(p)(\lambda_f(p)+\lambda_g(p))/2\big)}{p}\\
&\qquad\qquad\qquad\qquad+O\Big(\frac{1}{p^2} \Big)\bigg) +O\big((\log q)^{-10}\big)\bigg),
    \end{align*}
    by Lemmas \ref{lem:multiplicative}$(ii)$ and \ref{lem:gfest}$(iii)$. Again the product over $p \le c_0$ is $\ll 1$ and the second product over $y < p \le x$ is $\ll (\log \log q)^{O(1)}$ by Lemma \ref{lem:mediumprimeest}$(iii)$. For the third product over $c_0<p \le y$, we use \eqref{RSexp} to get that it is $\ll (\log y)^{-1/2}$. Hence,
\[
\mathbb E\Big( L^{f,g}(\tfrac12,X)\big|M_{\frac12}(X)\big|^2\Big)=O_\varepsilon\big( (\log q)^{-1/2+\varepsilon}\big),
\]    
and we obtain $(ii)$.
    \end{proof}

\bibliographystyle{amsalpha}

\end{document}